\documentclass[12pt]{article}
\usepackage{amsfonts}
\def\hybrid{\topmargin 0pt      \oddsidemargin 0pt
        \headheight 0pt \headsep 0pt
        \textwidth 16.5cm
        \textheight 23cm
        \hoffset=0.4cm
        \marginparwidth 0.0in
        \parskip 5pt plus 1pt   \jot = 1.5ex}
\catcode`\@=11
\def\marginnote#1{}

\newcount\hour
\newcount\minute
\newtoks\amorpm
\hour=\time\divide\hour by60
\minute=\time{\multiply\hour by60 \global\advance\minute by-\hour}
\edef\standardtime{{\ifnum\hour<12 \global\amorpm={am}%
        \else\global\amorpm={pm}\advance\hour by-12 \fi
        \ifnum\hour=0 \hour=12 \fi
      \number\hour:\ifnum\minute<10 0\fi\number\minute\the\amorpm}}
\edef\militarytime{\number\hour:\ifnum\minute<10 0\fi\number\minute}

\def\draftlabel#1{{\@bsphack\if@filesw {\let\thepage\relax
   \xdef\@gtempa{\write\@auxout{\string
      \newlabel{#1}{{\@currentlabel}{\thepage}}}}}\@gtempa
   \if@nobreak \ifvmode\nobreak\fi\fi\fi\@esphack}
        \gdef\@eqnlabel{#1}}
\def\@eqnlabel{}
\def\@vacuum{}
\def\draftmarginnote#1{\marginpar{\raggedright\scriptsize\tt#1}}

\def\draft{\oddsidemargin -0.1truein
        \def\@oddfoot{\sl preliminary draft \hfil
        \rm\thepage\hfil\sl\today\quad\militarytime}
        \let\@evenfoot\@oddfoot \overfullrule 3pt
        \let\label=\draftlabel
        \let\marginnote=\draftmarginnote
\def\@eqnnum{{\rm (\theequation)}
\rlap{\kern\marginparsep\tt\@eqnlabel}%
\global\let\@eqnlabel\@vacuum}  }

\newcommand{\RR}{{\mathbb{R}}}
\newcommand{\CC}{{\mathbb{C}}}

\newfont{\Bbbb}{msbm7 scaled 1\@ptsize00}
\newcommand{\zs}{\raise-1pt\hbox{$\mbox{\Bbbb Z}$}}

\font\sevenmsa=msam6 
\newfam\msafam
\textfont\msafam=\sevenmsa
\def\hexnumber@#1{\ifnum#1<10 \number#1\else
\ifnum#1=10 A\else\ifnum#1=11 B\else\ifnum#1=12 C\else
\ifnum#1=13 D\else\ifnum#1=14 E\else\ifnum#1=15 F\fi\fi\fi\fi\fi\fi\fi}
\def\msa@{\hexnumber@\msafam}
\def\llcorner{\delimiter"4\msa@78\msa@78 }
\def\lrcorner{\delimiter"5\msa@79\msa@79 }
\mathchardef\blacktriangleright="3\msa@49
\mathchardef\blacktriangleleft="3\msa@4A
\mathchardef\trianglerighteq="3\msa@44
\mathchardef\trianglelefteq="3\msa@45
\font\tenmsb=msbm10 scaled 1\@ptsize00
\newfam\msbfam
\textfont\msbfam=\tenmsb
\scriptfont\msbfam=\tenmsb
\def\msb@{\hexnumber@\msbfam}
\mathchardef\varkappa="0\msb@7B

\newdimen\linethick  \linethick=0.4pt
\newdimen\hboxitspace    \hboxitspace=5pt
\newdimen\vboxitspace    \vboxitspace=5pt

\def\fr#1{%
\beq\new
\vcenter{
\hrule height\linethick
           \hbox{\vrule width\linethick
                 \kern\hboxitspace
                 \vbox{\kern\vboxitspace
                       \hbox{$\begin{array}{c}\displaystyle#1
          \end{array}$}%
                       \kern\vboxitspace}%
                 \kern\hboxitspace
                 \vrule width\linethick}%
           \hrule height\linethick}%
\eeq}

\newdimen\Squaresize \Squaresize=14pt
\newdimen\Thickness \Thickness=0.5pt

\def\Square#1{\hbox{\vrule width \Thickness
   \vbox to \Squaresize{\hrule height \Thickness\vss
      \hbox to \Squaresize{\hss#1\hss}
   \vss\hrule height\Thickness}
\unskip\vrule width \Thickness}
\kern-\Thickness}

\def\Vsquare#1{\vbox{\Square{$#1$}}\kern-\Thickness}

\def\numberbysection{\@addtoreset{equation}{section}
        \def\theequation{\thesection.\arabic{equation}}}
\numberbysection

\renewcommand{\theequation}{\thesection.\arabic{equation}}
\newcommand{\l@qq}[2]{\addvspace{2em}
 \hbox to\textwidth{\hspace{1em}\bf #1 \dotfill #2}}

\newcounter{app}

\def\app{\setcounter{equation}{0}
\def\theequation{\Alph{app}.\arabic{equation}}\par
   \addvspace{10ex}
   \@afterindentfalse
  \secdef\@app\@dapp}
\newcommand\@app{\@startsection {app}{1}{-0.3ex}%
                             {-3.5ex \@plus -1ex \@minus -.2ex}%
                                   {2.3ex \@plus.2ex}%
                                   {\normalfont\Large\bf}}

\def\@dapp#1{%
{\parindent \z@ \raggedright \bf #1}\par\nobreak}
\def\l@app#1#2{\ifnum \c@tocdepth >\z@
    \addpenalty\@secpenalty
    \addvspace{1.0em \@plus\p@}%
    \setlength\@tempdima{1.5em}%
    \begingroup
      \parindent \z@ \rightskip \@pnumwidth
      \parfillskip -\@pnumwidth
      \leavevmode \bfseries
      \advance\leftskip\@tempdima
      \hskip -\leftskip
      #1\nobreak\hfil \nobreak\hb@xt@\@pnumwidth{\hss #2}\par
    \endgroup\fi}
\newcounter{sapp}[app]

\def\sapp{\def\theequation{\Alph{app}.\arabic{equation}}\par
   \@afterindentfalse
  \secdef\@sapp\@dsapp}
\newcommand\@sapp{\@startsection{sapp}{2}{\z@}%
                           {-3.25ex\@plus -1ex \@minus -.2ex}%
                           {1.5ex \@plus .2ex}%
                              {\normalfont\large\bfseries}}

\def\@dsapp#1{%
{\parindent \z@ \raggedright  \bf #1}\par\nobreak}

\newcommand{\l@sapp}{\@dottedtocline{2}{1.4em}{2.5em}}

\def\titlepage{\@restonecolfalse\if@twocolumn\@restonecoltrue\onecolumn
     \else \newpage \fi \thispagestyle{empty}\c@page\z@
        \def\thefootnote{\fnsymbol{footnote}} }

\def\endtitlepage{\if@restonecol\twocolumn \else  \fi
        \def\thefootnote{\arabic{footnote}}
        \setcounter{footnote}{0}}  
\relax

\hybrid

\newtoks\@stequation

\def\subequations{\refstepcounter{equation}%
  \edef\@savedequation{\the\c@equation}%
  \@stequation=\expandafter{\theequation}
  \edef\@savedtheequation{\the\@stequation}
  \edef\oldtheequation{\theequation}%
  \setcounter{equation}{0}%
  \def\theequation{\oldtheequation\alph{equation}}}

\def\endsubequations{%
  \setcounter{equation}{\@savedequation}%
  \@stequation=\expandafter{\@savedtheequation}%
  \edef\theequation{\the\@stequation}%
  \global\@ignoretrue}

\makeatletter
\newdimen\normalarrayskip            
\newdimen\minarrayskip               
\normalarrayskip\baselineskip
\minarrayskip\jot
\newif\ifold             \oldtrue            \def\new{\oldfalse}
\def\arraymode{\ifold\relax\else\displaystyle\fi}
\def\eqnumphantom{\phantom{(\theequation)}} 
\def\@arrayskip{\ifold\baselineskip\z@\lineskip\z@
     \else
     \baselineskip\minarrayskip\lineskip1\baselineskip\fi}

\def\@arrayclassz{\ifcase \@lastchclass \@acolampacol \or
\@ampacol \or \or \or \@addamp \or
   \@acolampacol \or \@firstampfalse \@acol \fi
\edef\@preamble{\@preamble
  \ifcase \@chnum
     \hfil$\relax\arraymode\@sharp$\hfil
     \or $\relax\arraymode\@sharp$\hfil
     \or \hfil$\relax\arraymode\@sharp$\fi}}

\def\@array[#1]#2{\setbox\@arstrutbox=\hbox{\vrule
     height\arraystretch \ht\strutbox
     depth\arraystretch \dp\strutbox
width\z@}\@mkpream{#2}\edef\@preamble{\halign \noexpand\@halignto
\bgroup \tabskip\z@ \@arstrut \@preamble \tabskip\z@ \cr}%
\let\@startpbox\@@startpbox \let\@endpbox\@@endpbox
  \if #1t\vtop \else \if#1b\vbox \else \vcenter \fi\fi
  \bgroup \let\par\relax
  \let\@sharp##\let\protect\relax
  \@arrayskip\@preamble}
\def\eqnarray{\stepcounter{equation}%
              \let\@currentlabel=\theequation
              \global\@eqnswtrue
              \global\@eqcnt\z@
              \tabskip\@centering              
              \let\\=\@eqncr
              $$%
            \halign to \displaywidth  \bgroup
             \eqnumphantom \@eqnsel
      \hskip\@centering                               
    $\displaystyle  \tabskip\z@ {##}$%
    &\global\@eqcnt\@ne \hskip 2\arraycolsep
         $ \displaystyle  \arraymode{##}$\hfil
    &\global\@eqcnt\tw@ \hskip 2\arraycolsep
         $\displaystyle\tabskip\z@{##}$\hfil
         \tabskip\@centering
    &{##}\tabskip\z@\cr}
\makeatother

\def\bea{\begin{eqnarray}}
\def\eea{\end{eqnarray}}

\def\bqa{\begin{eqnarray}}
\def\eqa{\end{eqnarray}}
\def\pr {\partial}

\def\beq{\begin{equation}}
\def\eeq{\end{equation}}
\def\be{\beq\new\begin{array}{c}}  
\def\ee{\end{array}\eeq}           
\def\bse{\begin{subequations}}                
\def\ese{\end{subequations}}

\def\stack#1#2{\raise0.7pt\hbox{$\mathrel{\mathop{#2}\limits^{#1}}$}}
\def\tr{\triangleright}                                
\def\tl{\triangleleft}                                 

\def\jo{\mathrel{\mkern-4mu}}
\def\sem{\mathsurround=0pt
\mathrel{\raise1.4pt\hbox{$\scriptscriptstyle>$}}\jo\mathrel\tl}
\def\mes{\mathsurround=0pt
{\mathrel\tr\jo\mathrel{\raise1.4pt\hbox{$\scriptscriptstyle <$}}}}

\def\]{]\raise-2pt\hbox{$_\ast$}}
\def\op#1{\raise-6pt\hbox{$\stackrel{\displaystyle\oplus }{\scriptstyle
#1}$}\;}
\def\oop#1#2{\raise-6pt\hbox{$\stackrel{#2}{\stackrel{\displaystyle\oplus}
{\scriptstyle #1}}$\;}}

\def\la{\lambda}

\def\<{\langle}
\def\>{\rangle}

\def\wt{\widetilde}
\def\wh{\widehat}
\def\t{\theta}
\def\frak{\mathfrak}
\def\g{\gamma}

\def\gfr{\frak{g}}

\def\N{{\scriptscriptstyle N}}
\def\ts#1#2{{\textstyle\frac{#1}{#2}}}

\begin{document}
\hfill ITEP-TH-80/02
\begin{center}
\thispagestyle{empty}
\phantom.
\bigskip\bigskip\bigskip\bigskip
{\Large\bf On a class of  integrable systems connected with
$GL(N,\RR)$\footnote{Extended talk by the third author at the
6 th International Workshop: "Conformal Field Theory and Integrable
Models", September 15-21, 2002, Chernogolovka, Russia.}}

\bigskip

\bigskip \bigskip
{\large A. Gerasimov\footnote{E-mail: gerasimov@vitep5.itep.ru},
 S. Kharchev\footnote{E-mail:  kharchev@gate.itep.ru}, D.
Lebedev\footnote{E-mail:  lebedev@gate.itep.ru}}\\ \medskip
{\it Institute of Theoretical \& Experimental Physics\\ 117259 Moscow,
Russia}\\
\end{center}

\vspace{4cm}

\begin{abstract}
\noindent In this paper we define a new class of the quantum
integrable systems associated with the quantization of the
cotangent bundle $T^*(GL(N))$ to the Lie algebra $\frak{gl}_N$.
The construction is based on the Gelfand-Zetlin maximal commuting
subalgebra in $U(\frak{gl}_N )$. We discuss the
 connection with the other known integrable systems based on
$T^*GL(N)$. The construction of the spectral tower associated with
the proposed integrable theory  is given. This spectral tower
appears as a generalization of the standard spectral curve for
integrable system.

\end{abstract}

\clearpage \newpage



\normalsize
\section{Introduction}

Let us start with an informal general definition of the quantum
integrable system. General quantum dynamical system may be
described by the quadruple $(A,H,\pi,\mathcal{H})$ where $A$ is a
$C^*$-algebra, Hamiltonian $H$ is an element of $A$ defining the
evolution  of the quantum system, $\mathcal{H}$ is a Hilbert space
of states of the system and $\pi$ is a faithful representation of
$A$ by the unitary operators acting in $\mathcal{H}$.

 The integrable structure  on the dynamical system $(A,H,\pi,\mathcal{H})$
is a choice of the maximal commutative subalgebra $I\subset A$
such that $H\in I$. In the non-degenerated case the maximality is
equivalent to the conditions on the dimension of the subalgebra
$I$ to be equal to the half of the dimension of $A$.  When the
algebra $A$ is included in the flat family of the algebras
$A_{\hbar}$  such that  $A_0$ is a commutative algebra the last
condition means ${\rm dim(Spec}(I))=\frac{1}{2}{\rm
dim(Spec}(A_0))$. The subalgebra $I$ is the algebra of integrable
Hamiltonians. Up to some details the Hilbert space $\mathcal{H}$
may be realized as a (sub)space of functions on some configuration
space $X$: $\mathcal{H}\sim L^2(X)$. Thus we have another natural
maximal commutative subalgebra $I_0\subset A$ - the algebra of
functions on $X$. In interesting cases the subalgebras $I$ and
$I_0$ are different and in particular $H$ is not in $I_0$. The
solution of the quantum integrable system usually means an
explicit construction of the common eigenvalues of the commutative
subalgebra $I$ acting in $\mathcal{H}$. The maximality condition
for $I$ implies that the eigenspaces are one-dimensional for the
generic point of ${\rm Spec}(I)$. It is clear that the explicit
solution is given by some kind of Fourier transformation from the
subalgebra $I$ to the subalgebra $I_0$.  The choice of the
particular Hamiltonian $H\in I $ will not be important in the
following discussion and we will describe the integrable system as
a quadruple: $(A,I,\pi,\mathcal{H})$.

In this note we describe a class of the quantum integrable systems
associated with Lie groups. Our main example will be the
integrable system associated with the quantization of the
cotangent bundle $T^*GL(N)$ to the Lie group $GL(N)$. We define
the maximal commuting subalgebra of the quantized algebra of
functions on $T^*GL(N)$. Also some examples of the quantum
reductions respecting the defined integrable structure will be
shortly discussed.

Starting with the pioneering works by Olshanetsky and Perelomov
\cite{OP1,OP2,OP3} the constructions of the integrable systems
connected with the Lie group  $G$ corresponding to the Lie algebra
$\frak{g}$ follow more or less the same  procedure. Considering
the quantized algebra of functions on $T^*G$ one has an obvious
choice of the commutative subalgebra - the center of the universal
envelope $\mathcal{Z}\subset U(\frak{g})$. However the functional
dimension of this commutative subalgebra is too small to define
the integrable structure on $T^*G$. According to the standard
reasoning to obtain the integrable system  one should consider the
reduction under the appropriate (left/right) actions of the group
$G$ so that the center $\mathcal{Z}$
 descends under reduction to the commutative subalgerba which is
 large  enough to provide the integrable structure for the reduced
 system. Conceptually this is not quite
 satisfactory. The integrable structure appears only at the final step and
 the dimensionality condition provides strong constraint
on the possible reductions.

To cure this deficiency one would like to define the integrable
structure on $T^*G$ before the reduction. One such construction is
known for a long time.  In the classical approximation (i.e. the
maximal commuting family with respect to the Poisson brackets) it
was proposed by Mishchenko \cite{M}(see also \cite{Man},
\cite {MF}). Later basically the same constructions appears in
\cite{FR,FT,AHH} as a particular case of the constructions of
integrable theories connected with loop groups. The center of loop
algebra was used to define the integrable structure on the general
orbits of the finite-dimensional groups.

In this note we follow another route and define the different
non-trivial maximal commutative subalgebra of the algebra of
quantized functions on $T^*GL(N)$. The main ingredient is
Gelfand-Zetlin construction of the maximal commutative subalgebra
in $U(\frak{gl}_N)$ \cite{GelZ}, \cite{Cher1,Cher2}. Obviously
thus defined  integrable systems are closely related with the
representation theory of the groups. In particular matrix elements
in the Gelfand-Zetlin representation \cite{GelG} provides the
explicit description of the quantum wave functions for the
corresponding integrable system.

Non-trivial integrable systems  usually associated with
interesting spectral geometry. The dynamics is linearized on the
Jacobi variety of the corresponding spectral curve. We will show
that the proposed integrable structures on $T^*GL(N)$ (and general
coadjoint orbits of the group) give rise to the interesting
generalization of the spectral geometry- {\em spectral tower}. It
seems remarkable that we could associate a spectral geometry with
the group itself (or better to say with $T^*GL(N)$) and not with
its special reductions.

Let us stress, that the subalgebra which is used in this note to
define integrable structure, is not contained in the center of
the loop algebra in any obvious way. On the contrary one would
like to combine the proposed construction with Mischenko-Fomenko
approach and give a unified framework for the description of
various integrable structures on cotangent bundles to  (loop)
groups. We hope to return to these questions elsewhere.

The rest of the paper consists of two parts. In  section  2 we
start with the construction of the Poisson commuting families of
functions on
 $T^*GL(N)$ which  provide   the classical integrable structures.
Then we discuss the quantization of this algebra. At the end of
this  section some examples of the reductions consistent with the
given integrable structure are considered. This leads in
particular to the integrable systems associated with the coadjoint
orbits of $GL(N)$ (which implicitly appeared in \cite{GKL}).

In section 3 we construct an analog of the "spectral curve"
connected with the integrable system introduced in section 2. It
turns out that the resulting geometry most adequately described by
the notion of {\em spectral tower}. It seems closely connected
with the Bott tower and Bott-Samelson manifolds in  differential
geometry and representations theory \cite{BS,BT,GK}.

\section{Integrable structures on $T^*GL(N)$}

In this section we introduce the family of the non-trivial
integrable structures based on the quantum deformation
$\mathcal{A}$ of the algebra of functions on $T^*G$ for $G=GL(N)$.
As it was discussed in the introduction to define the integrable
structure one should chose the maximal commutative subalgebra
$I\subset \mathcal{A}$. We  start with the description of the
corresponding classical integrable structure which is defined by
the choice of the two maximal commutative subalgebras in the
Poisson algebra of functions on $T^*G$.

Consider right action of the group $G$ on its  cotangent bundle
$T^*G$. Using this action one could  identify $T^*G$  with the
product of the group $G$ and the dual $\frak{g}^*$ to the
corresponding Lie algebra $\frak{g}$. Using the non-degenerate
pairing on the Lie algebra $\gfr=\frak{gl}_N$ we may also identify
$\gfr^*\cong \gfr$. Thus the points in $T^*G$ are described by the
coordinates $(u,g)\in \gfr \times G$. If we would use the left
action to identify the cotangent bundle with the product  of
the Lie algebra and the Lie group we will get another
parameterizations in terms of the coordinates
$(\widetilde{u},\widetilde{g}) \in \gfr \times G$. The connection
between the different coordinates is given by
$\widetilde{u}=gug^{-1},\widetilde{g}=g$.

The cotangent bundle is canonically supplied with the symplectic
structure which has the following form:
\bqa
\omega=<du,dgg^{-1}>+<u,dgg^{-1}\wedge dgg^{-1}>\,.
\eqa
Corresponding Poisson structure is defined by the basic relations:
 \bqa \label{Poisson00}
 \{u^{ij},u^{kl}\}&=&\delta^{jk}u^{il}-\delta^{li}u^{kj}\,,\\
\label{Poisson01}
\{\widetilde{u}^{ij},\widetilde{u}^{kl}\}&=&
\delta^{jk}\widetilde{u}^{il}-\delta^{li}\widetilde{u}^{kj}\,,\\
\label{Poisson1} \{u^{ij},\widetilde{u}^{kl}\}&=&0\,,\\
 \{u^{ij};g_{kl}\}&=& \delta_{il}g_{kj}\,,\\
\{\widetilde{u}^{ij};g_{kl}\}&=& \delta_{jk}g_{il}\,,\\
 \label{Poisson2} \{g_{ij},g_{kl}\}&=&0\,.
 \eqa
Thus the Poisson brackets (\ref{Poisson00}) and (\ref{Poisson01})
reproduce the underlying Lie algebra structure. Note that
 $u^{ij}$ and $\widetilde{u}^{kl}$ play the role of  the momentum
maps for the left and right actions of $G$ on $T^*G$ and thus
(\ref{Poisson1}) is pretty obvious.

As a "coordinate" commutative subalgebra $I_0$ we chose the
evident maximal Poisson commutative subalgebra generated by
$g_{ij},\;i,j=1,\ldots, N$. So the corresponding quantum dynamical
system will be defined in terms of the representation in the
Hilbert space $L^2(G)$ of the square integrable functions on the
group  $G$. To define the non-trivial integrable theory we should
construct another maximal commutative subalgebra $I$ which is in
general position with $I_0$. Rather obviously the following set of
functions \bqa I^u_{ij}&=&\sum_mg^{-1}_{im}u_{mj}\eqa generates
the maximal commutative subalgebra. However this subalgebra is not
very interesting from the point of view of the quantum integrable
theories due to its simple connection with $I_0$.

 Taking into account the results of \cite{GelZ,GelG}
(see also \cite{GKL}   for applications to integrable systems) the
natural guess for the non-trivial Poisson commutative subalgebra
is the following. Consider the  polynomials: \bqa \label{Gelf}
I^L_k(\la)&=&\sum_nI^L_{kn}\la^n=\det(\lambda-u)_k\,, \ \ \ \ \
k=1, \ldots,N-1\,,\\
I^R_k(\la)&=&\sum_nI^R_{kn}\la^n=\det(\lambda-\widetilde{u} )_k\,,
\ \ \ \ \ k=1, \ldots,N-1\,, \\
I_N(\la)&=&\sum_nI_{Nn}\la^n=\det(\lambda-u), \eqa where ${\rm
det}(A)_k$ is $k$-th principal  minor of the matrix $A$. It is not
difficult to see that the algebra generated by the coefficients of
these polynomials is commutative.  Note that $I_N(u,\la)$ is
invariant under conjugation of $u$ by any group element and thus
in particular under the Hamiltonian flows generated by the
elements $I_{Nl}$ for any $l$. Thus the family $I_{Nl}$ is Poisson
commutative. By the same reason we have
$\{I_{N,l},I^L_{N-1,k}\}=0$ for any $k=1,\ldots,N-1$. It is easy
to see from (\ref{Poisson00}) that Poisson brackets of $u^{ij}$
for $i,j=1,\ldots,N-1$ reproduce the relations of the Lie algebra
$\frak{gl}(N-1)$. Thus  by the same reason  $I^L_{N-1,k}$ commute
with any $u^{ij}$ for $i,j=1,\ldots,N-1$. Hence  we have the
commuting family generated by $I_{Nk}$,$I^L_{N-1,l}$. Applying the
induction over the rank we derive the commutativity of the algebra
generated by $I_{Nk}$,$I^L_{nl}$ $n=1,\cdots ,N-1$. Taking into
account (\ref{Poisson00}),(\ref{Poisson01}), (\ref{Poisson1}) we
have proved that the algebra generated by the coefficients of the
polynomials (\ref{Gelf}) is a commutative algebra. The number of
the  independent coefficients of these polynomials is $N^2$ ( half
of the dimension of
 $T^*G$)  and  thus the subalgebra generated by the coefficients of
 the polynomials satisfies the dimensionality condition.
 Together with the choice of the "coordinate"
subalgebra $I_0$ this provides us with the classical integrable
structure on $T^*G$.

To construct the maximal commutative subalgebras  one  could use
the determinants of the other families of submatrices. For
instance instead of the left-upper corner minors we may use the
left-lower corner minors of the matrices $\la-u$ and
$\la-\wt u$. For the description of the various
reductions it will useful in the following to characterize the
possible choices of the family of the submatrices in terms of the
invariance properties. Thus the principal minors of the matrix $M$
are invariant under the transformations $M\rightarrow N_+MN_-$
with $N_{\pm}$ being strictly lower/upper triangle matrices and
the left-lower corner minors which are invariant with respect to the
transformations $M\rightarrow N_+MN_+$.

Let us now consider the corresponding quantum integrable theories.
We define the quantization of the Poisson commutation relation
together with its faithful representation in $\mathcal{H}\cong
L^2(G)$ as follows. The quantum counterparts $u_q$ and
$\widetilde{u}_q$   of $u$ and $\widetilde{u}$ act by the left and
right differentiation: \bqa u_q^{ij} = \nabla^{ij}_L\equiv\sum_k
g_{ki}\frac{\pr}{\pr g_{kj}}\,,\\ \widetilde{u}_q^{ij} =
\nabla^{ij}_R\equiv -\sum_k g_{jk}\frac{\pr}{\pr g_{ik}} \eqa and
$g_{ij}$ act by multiplication. This operators obviously satisfy
the relations (\ref{Poisson00})-(\ref{Poisson1}) with commutators
instead of the Poisson brackets.

 To quantize the integrable system it is necessary to
describe the quantum counterpart of the maximal commutative
subalgebra $I$ (the quantization of the subalgebra $I_0$ is
obvious in our case).  Given any construction of the center
$\mathcal{Z}$ of the universal envelope
 $U(\gfr)$ one could produce the maximal commutative subalgebra
following the original proposal \cite{GelZ} (the similar classical
arguments were described above). The construction of the
generators $\mathcal{Z}$ is obviously not unique. For instance one
could use  the symmetrization map $s: x_1x_2\cdots x_n \rightarrow
\frac{1}{n!} \sum_{\sigma \in S_n}x_{\sigma_1}x_{\sigma_2} \cdots
x_{\sigma_n}$ from the invariants of the adjoint action of $\gfr$
on the symmetric powers $S^{\bullet}(\gfr)$ to the invariants of
the adjoint action on $U(\gfr)$. Also, one could use the Duflo
isomorphism \cite{Du}. Taking into account that classical
construction was given in terms of the characteristic polynomials
of various submatrices it would be more appropriate to use another
construction based on the notion of the quantum determinant
\cite{KS}. It is well known that the coefficients of the following
polynomial:
\be
{\rm qdet}(\la-u) \\ =\,\sum_{p\in P_\N}{\rm sign}\,p\;
(\la-\rho^{(\N)}_1-u)_{p(1),1}\ldots(\la-\rho^{(\N)}_k-u)_{p(k),k}
\ldots (\la-\rho^{(\N)}_\N-u)_{p(\N),\N}\,,
\hspace{-0.5cm}
\ee
where $\rho^{(\N)}_n=\ts{1}{2}(N-2n+1),\,(n=1,\ldots,N)$ and
summation is given over elements of the permutation group
$P_\N$, generate the center $\mathcal{Z}\subset U(\mathfrak{gl}(N))$.
Now using the iterative construction described in the classical
case we obtain the quantum maximal  commutative subalgebra
generated by the coefficients of the polynomials:
\bqa
\label{qGelf}
 I^L_k(\la)&=&\sum_nI^L_{kn}\la^n={\rm qdet}(\lambda-u_q)_k\,,
\ \ \ \ \ \ \ k=1, \ldots, (N-1), \\
I^R_k(\la)&=&\sum_nI^R_{kn}\la^n={\rm
qdet}(\lambda-\widetilde{u}_q )_k\,, \ \ \ \ \ \ \ k=1, \ldots,
(N-1) ,\\ I_N(\la)&=&\sum_nI_{Nn}\la^n={\rm qdet}(\lambda-u_q)\,.
\eqa

The explicit solution of the quantum integrable systems based on
the Gelfand-Zetlin construction has a clear meaning in the
representation theory. Consider the quantum wave function
$\Psi_{C}(g)$ which is a solution of the  equation eigen-function
 equation:
\bqa\label{eigen}
\wh I^{\;L,R,N}_{nm}\Psi_{C}=C^{L,R,N}_{nm}\Psi_C\;.
\eqa
We use the subalgebra $I_0$ generated by $g_{ij}$ to define the
"coordinates". Thus the wave function is supposed to be an element
of the Hilbert space $\mathcal{H}\simeq L^2(G)$ and
quantum Hamiltonians are realized  $\wh I^{\,L,R,N}_{nm}$ in terms
of the differential operators acting on $L^2(G)$ with the
eigenvalues  $C^{L,R,N}_{nm}$. Being written in this form, the
equations (\ref{eigen}) represent a set of the differential
equations on $\Psi_C(g)$. The explicit solution may be constructed
in terms of the representation theory. Let us introduce the new
variables $\gamma_{nk}$ and $\wt\gamma_{nk}$ as
follows:
\bqa
 I^L_k(\la)&=&\prod_n (\la-\gamma_{kn})\,,\\
 I^R_k(\la)&=&\prod_n (\la-\wt\gamma_{kn})\,,\\
 I^N(\la)&=& \prod_n (\la-\gamma_{Nn})\,.\eqa
One could use the new variables $\gamma,\widetilde{\gamma}$ as
"coordinates" to realize the representation of $GL(N)$ in terms of
the difference operators acting on the functions of the variables
$\g$ or $\widetilde{\g}$. This construction was described in
details in \cite{GKL}. Consider the matrix element of the group
element $g$ acting in the representation $T$:\bqa
<\gamma|T(g)|\widetilde{\gamma}>\eqa Clearly this function of $g$
satisfies the equations (\ref{eigen}) with: \bqa
 C^L_k(\la)&=&\prod_n (\la-\gamma_{kn})\\
 C^R_k(\la)&=&\prod_n (\la-\widetilde{\gamma_{kn}})\\
 C^N(\la)&=& \prod_n (\la-\gamma_{Nn})\eqa
Thus we have a clear procedure to get a representation for the
wave function of the quantum integrable system on $T^*G$.

Now we briefly describe the various reductions
 compatible with the integrable structures on $T^*G$
 defined above. We start with the general description
 on the quantum reduction. Let us given a symplectic manifold $M$
 and let $\mathcal{A}_{\hbar}(M)$ be a quantization of the algebra
 of functions on $M$ supplied with the faithful representation
$\pi: \mathcal{A}_{\hbar}(M)
 \rightarrow {\rm End}\,(\mathcal{H})$. Suppose we also fix some integrable
 structure defined by the maximal commutative subalgebra $I \in
 \mathcal{A}_{\hbar}(M)$. Consider the action of the Lie algebra
 $\frak{b}$ on $\mathcal{A}_{\hbar}(M)$ by the inner
 derivatives. Classically it corresponds to the Hamiltonian
 action of $\frak{b}$ on $M$. Thus given the integrable structure
$(\mathcal{A}_{\hbar}(M),I,\pi,\mathcal{H})$ we define the reduced
integrable structure as
$(\mathcal{A}^{Red}_{\hbar}(M),I^{Red},\pi^{Red},\mathcal{H}^{Red})$
where \bqa
 \mathcal{A}^{Red}_{\hbar}(M)&=&\{ a\in \mathcal{A}^{Red}_{\hbar}(M)
|\,[a,b]=0; a\sim a+b, \forall b\in \frak{b}\},\\
 \mathcal{H}^{Red}&=&\{v\in \mathcal{H}|\,bv=0,\forall b\in \frak{b}\}.\eqa
The reduced commutative subalgebra $I^{Red}$ is constructed using
the explicit choice of the generators $(I_1,\ldots I_M)$ of $I$.
It is generated by the $\frak{b}$-invariant subset of the
generators $(I_{r_1},\ldots ,I_{r_M})$.  When $I^{Red}$ appears
 to be maximal commutative subalgebra we  say that the
reduction is compatible with the integrable structure on $I$.
Obviously this definition of the compatibility of the reduction
and the integrable structure is too restrictive. However in the
following we will consider only the examples that satisfy this
condition.

Now consider the examples of the quantum reduction of the
integrable structure on $T^*G$.  The most simple case is the
(quantum) Hamiltonian reduction under the left (or right) action
of $G$ on $T^*G$.  The reduced phase space is a coadjoint orbit
$\mathcal{O}$ of $G$. The polynomials $I_N(\la),I^L_k(\la)$
($I_N(\la),I^R_k(\la)$) are invariant with  respect to the left
(right) action of the group and
 thus define the integrable structure on the coadjoint orbit.
 This set of commuting elements gives a maximal
commuting subalgebra and  provides the integrable structure on the
algebra of quantized functions on the orbit. Actually this
integrable structure was used in \cite{GKL} to define
Gelfand-Zetlin parameterization of the coadjoint orbit. One could
consider the reduction under the lift to $T^*G$ of the adjoint
action of $G$. In this case the classical reduced phase space is a
 factor $T^*H/W$ of the cotangent bundle to the Cartan torus $H$
of $G$ under the action of the Weyl group $W$. Here the integrable
structure on the reduced theory is given by the commuting
subalgebra $I^N_n$.

Another interesting class of the integrable systems is obtained
using  the reductions under the combination of the left and right
actions of various subgroups of $G$ (which leads to the notion of
$R$-groups \cite{RS}).  The particular case is the reduction under
the action of the algebra of the low triangle matrices from the
left and upper-triangle matrices from the right. If the
corresponding orbits are zero-dimensional we get the open Toda
chain phase space with the integrable structure defined by $I_N$.
This is the original construction of the open Toda chain in the
group theory approach due to Kostant \cite{Ko} (see also
\cite{STS-toda}). The connection between the reduced and unreduced
quantum integrable structures was used in \cite{GKL} to construct
a representation of the wave function for open Toda chain
\cite{KL2} in terms of the representation theory. Let us also
mention the interesting intermediate case  when the orbit of the
left action is a point and the orbit of the right action of the
triangle matrix is generic.  The corresponding space is a twisted
cotangent bundle to the Borel subgroup $T^*_fB$. To define the
integrable structure on the reduced space  it is convenient to
use another choice of the integrable structure on $T^*G$
constructed with the help of the left-lower minors of $(u-\la)$
and $(\widetilde{u}-\la)$. This case was considered previously in
\cite{D}.

At the end of this section let us mention another non-trivial
choice of the maximal commuting family.  This family was
introduced by Mishchenko in \cite{M} (see also \cite{MF}) and is
based on the  commutative algebra  generated by the coefficients
$I^A_{nm}$ of the polynomial in two formal variables $\la$ and
$\mu$: \bqa \label{MF}
 I^A(\mu,\la)=\sum_{nm}I^A_{nm}\la^n\mu^m=\det(u-\mu A-\la)\,,
\eqa where $A$ is some fixed matrix. This construction provides a
 different maximal commutative subalgebra. It has a natural
  interpretation in terms of the loop algebra $LG$
associated with the group $G$. The functional (\ref{MF}) turns out
to be the restriction of the generator of the center of the loop
algebra:\bqa
 I^A(\mu,\la)=\sum_{nm}I^A_{nm}\la^n\mu^m=\det(u(\mu) -\la)\eqa
on the Poisson manifold obtained by the reduction of the cotangent
bundle $T^*LG$ under the combined left/right action of some
appropriate subgroups of $LG$. There are indications that the
construction proposed in this note  being properly generalized to
 loop group case  provides a unified description of the various
commutative family including the Gelfand-Zetlin and Mischenko
constructions. We are going to consider this generalization in the
future publication.

\section{Spectral tower associated with  $T^*GL(N)$}

 In this section we show that it is possible to associate the
 family of the spectral curves with the classical integrable structure
 constructed in the previous section. To
simplify the presentation we deal not with the integrable
structure constructed on $T^*G$ but consider the reduction under
(left) action of the group $G$ over the generic orbit. The
geometry of the family of the spectral curves associated with the
coadjoint orbits  has the new interesting properties and most
adequately described by the notion of the {\em spectral tower}
which generalizes the notion of the usual spectral curve of the
integrable systems (see for example \cite{DMN}, \cite{Mum}). As
was discussed previously  the integrable structure on the generic
orbit $\mathcal{O}\subset \frak{gl}(N)$ descends from the
integrable structure on $T^*GL(N)$ and is described by the
polynomials: \bqa\label{orb} I^L_k(\la)=\det(\lambda-u)_k\,, \eqa
where $ k=1, \ldots, N$.  It turns out that the each polynomial
$I^L_k(\la)$  defines the corresponding spectral curve $\Sigma_k$.
However the whole integrable structure does not naturally split
into the factors corresponding to the spectral curves associated
with each polynomial. The point is that fixing spectral invariants
of the left-upper $(k-1)\times (k-1)$ submatrix of $u$ we also
partially fix the representative in the conjugation class of the
left-upper $k\times k$ submatrix of $u$. This subtlety is a
manifestation of the non-trivial choice of the zero section in the
family of the Jacobians of the spectral curves and results in the
iterated structure of the spectral tower.

To uncover the spectral geometry behind the integrable theory
associated with the coadjoint orbit  we start with the symplectic
form in  the Gelfand-Zetlin canonical variables and then represent
the symplectic structure in the form that has  the obvious
interpretation in terms of the spectral curves. First we introduce
the  Gelfand-Zetlin variables $(\g_{nj},\t_{nj})$ for $GL(N)$. The
details of its construction may be found in \cite{GKL} (see also
\cite{AFS,S} for the case of the compact groups). Also to make
contact with \cite{GKL} we change the notations $A_n(\la)\equiv
I^L_n(\la)$. Define the variables $\g_{nk}$ as the roots of the
polynomials $A_n(\la)$: \bqa
A_n(\la)&=&\prod_{j=1}^n(\la-\g_{nj})\,,\ \ \ \ \ n=1,\ldots,N\,,
\eqa and $\t_{nj}$ are the conjugated variables so that the
canonical Kirilov-Kostant symplectic form on the orbit is given:
\bqa \omega=\sum_{n,j} \delta \t_{nj}\wedge \delta \g_{nj}\,. \eqa
It is useful  to introduce the following polynomials $C_n(\la)$:
\bqa C_n(\la)=-\sum_{j=1}^n
\prod\limits_{r=1}^{n-1}(\gamma_{nj}-\gamma_{n-1,r}) \prod_{s\neq
j}\frac{\la-\gamma_{ns}}{\gamma_{nj}-\gamma_{ns}}\;
 e^{\theta_{nj}}\,.
\eqa Thus, in particular we have \bqa
C_n(\g_{nj})=-A_{n-1}(\g_{nj})e^{\t_{nj}}\,. \eqa Simple
calculations lead to the following expression for the symplectic
structure in terms  of $A_n(\lambda)$ and $C_n(\lambda)$: \be
\label{CAA} \omega=\\ \frac{1}{2\pi i} \sum_n \int_{\Gamma^C_n}
\delta \log C_n(\lambda) \wedge \delta \log A_n(\la) -
\frac{1}{2\pi i} \sum_n \int_{\Gamma^A_{n-1}} \delta \log
A_{n-1}(\lambda) \wedge \delta \log A_n(\la). \ee Here the
contour $\Gamma^C_n$ ($\Gamma^A_{n-1}$) is a sum of the small
contours around zeros of the polynomial $C_n(\lambda)$
($A_{n-1}(\lambda)$).

Let us define the new variables $h_k$ as the coefficients of
$A(\la)$:
\bqa\label{an}
 A_n(\la)=\sum_k \la^{n-k} h_{nk}\,.\eqa
Note that these variables are invariant with respect to the action
of the permutation group $P_n$ acting on the roots of the
polynomial $A_n(\lambda)$. Thus  although the variables $\g_{nk}$
are defined on the finite cover of  (the open part of) the orbit
the variables $h_{nk}$ parameterize the orbit itself.
 It is easy to verify that the symplectic form has the canonical form:
 \bqa \omega&=&\sum_{n,j} \delta \tau_{nj}\wedge \delta h_{nj}
 \eqa
in the coordinates $(h_{nk},\tau_{nk})$:
\bqa
 h_{nk}&=&\frac{1}{2\pi
 i}\oint_{\infty}A_n(\lambda)\frac{d\la}{\la^{n-k+1}}\\
\tau_{nk}&=&\sum_{i=1}^{n-1}\int ^{e_{ni}}_{\la_0}
\frac{\la^{n-k}} {A_n(\la)}\,d\la - \sum_{i=1}^{n-1}\int
^{\g_{n-1,i}}_{\la_0} \frac{\la^{n-k}} {A_n(\la)}\,d\la\;.
\eqa
Here
$\la_0$ is an arbitrary point on the $\la$-plane and $e_{nj}$ are
zeros of the polynomials $C_n(\la)$.

Now we could easily recover  the spectral curve behind the
integrable system on the coadjoint orbit. Let $\Sigma_n$ be a
complex non-compact curve obtained from the projective plane $\CC
P^1$ by deleting the points $(\g_{ni})$,$i=1, \ldots ,n$. The
spectral curve $\Sigma$ of the integrable system will be  a
disjoint union of the complex curves $\Sigma_n$.  The basis of the
holomorphic one-differentials on $\Sigma_n$ that is naturally
connected with the parameterization of the curve by the
coefficients $h_{nj}$  of the polynomials (\ref{an}), is given by
\bqa \label{bases}
\Omega^{(k)}_n=\frac{\la^{k-1}}{A_n(\la)}\,d\la\,. \eqa The
corresponding (generalized) Jacobi variety has a very simple
description in this case: $J(\Sigma)=(\CC^*)^n$ (it is obvious in
another basis $\omega_{nj}$ normalized by the condition ${\rm
Res}_{\la=\g_{nk}}\omega_{ni}=2\pi i\delta_{kj}$). Let ${\cal
B}_n=(\CC P^1)^n/\Delta_n$ be a configuration space of the
disjoint points $\g_{n1},\ldots ,\g_{nn}$ in $\CC P^1$.  The
family of the spectral curves $\Sigma=\bigsqcup \Sigma_n$
parameterized by $\{\g_{nj}\}\in \prod_{n=1}^N{\cal B}_n$ gives
rise to the family of the
 Jacobians $J(\Sigma)=\prod_{n=1}^N J(\Sigma_n)$,
$J(\Sigma_n)=(\CC^*)^n$.
 The total space of the family  is the phase space of our integrable system.
 The variables $h_{nk}$ are naturally interpreted as
 the "action"-type variables which parameterize the base of the
 fibration i.e. the moduli space of the curve $\Sigma$
 and the variables
 $\tau_{nk}$ are the image of the divisors ${\cal D}_n=\sum_{i=1}^n
 (e_{ni})-\sum_{j=1}^{n-1}(\g_{n-1,j})-(\la_0)$ under the Abel
 map:
\bqa
 \la \rightarrow \int^{\la}_{\la_0}\Omega_n^{(k)}\,.
\eqa
Note that the natural zero section $\rho:{\cal B}_{n-1}\rightarrow
J(\Sigma_n)$ of the family of Jacobians $J(\Sigma_n)$ is given the
image of the divisor  $\sum_{j=1}^{n-1}(\g_{n-1,j})+(\la_0)$
\bqa
\rho: \{\g_{n-1,j}\}\rightarrow \{\sum_j^{n-1}
 \int^{\g_{n-1j}}\Omega_n^{(k)}\}\,.
\eqa Thus although the family structure is factorized on the
product of the families connected with each $\Sigma_n$ the section
$\rho$ does not. The family with this  section has the structure
of the iterated fibration: \bqa\label{Tower} M_n\rightarrow
M_{n-1} \rightarrow M_{n-2} \rightarrow \ldots \rightarrow M_1\;,
\eqa where the fiber of $M_k\rightarrow M_{k-1}$ is a total space
of the family of the Jacobians associated with $\Sigma_k$. It
seems reasonable to coin the name {\em spectral tower} for this
structure. It is the non-trivial structure of the spectral tower
that is responsible for the appearance of the second term in the
expression for the symplectic structure (\ref{CAA}).

Let us note finally that the described  structure is in a some
sense  similar to the construction of the splitting of a vector
bundle over the differential manifold \cite{BS,BT}. For the close
constructions in representation theory see \cite{GK}.

\section*{Acknowledgments}
We would like to thank the Max Planck Institut f\"ur Mathematik
(Bonn, Germany), where the paper was partially done, for
stimulating atmosphere. The research was partly supported by
grants RFBR 00-02-16530 (A.Gerasimov); INTAS OPEN 00-00055; RFBR
01-01-00539 (S.Kharchev); INTAS OPEN 00-00055; RFBR 00-02-16530
(D.Lebedev) and by grant 00-15-96557 for Support of Scientific
Schools.

\end{document}